\documentclass[final,onefignum,onetabnum]{siamart171218}

\usepackage{lipsum}
\usepackage{amsfonts}
\usepackage{graphicx}

\ifpdf
  \DeclareGraphicsExtensions{.eps,.pdf,.png,.jpg}
\else
  \DeclareGraphicsExtensions{.eps}
\fi

\usepackage[noend]{algpseudocode}

\algnewcommand{\LineComment}[1]{\vspace{.0625in}\State \textcolor{gray}{\texttt{\# #1}}}
\algrenewcomment[1]{\hfill\textcolor{gray}{\texttt{\# #1}}}

\newsiamremark{remark}{Remark}
\newsiamremark{hypothesis}{Hypothesis}
\crefname{hypothesis}{Hypothesis}{Hypotheses}
\newsiamthm{claim}{Claim}

\headers{Concurrent Multi-Parameter Learning}{B. Pachev, J. P. Whitehead, and S. A. McQuarrie}

\title{Concurrent multi-parameter learning demonstrated on the Kuramoto-Sivashinsky equation}

\author{Benjamin~Pachev\thanks{Oden Institute for Computational Engineering and Sciences, University of Texas at Austin, Austin, TX (\email{benjaminpachev@gmail.com}, \email{shanemcq@utexas.edu}).}
\and Jared~P.~Whitehead\thanks{Department of Mathematics, Brigham Young University, Provo, UT \\
  (\email{whitehead@mathematics.byu.edu}).}
\and Shane~A.~McQuarrie\footnotemark[1]}

\usepackage{amsopn}

\ifpdf
\hypersetup{
  pdftitle={Concurrent multi-parameter learning demonstrated on the Kuramoto-Sivashinsky equation},
  pdfauthor={B. Pachev, J. P. Whitehead, S. A. McQuarrie}
}
\fi

\begin{document}

\maketitle

\begin{abstract}
We develop an algorithm based on the nudging data assimilation scheme for the concurrent (on-the-fly) estimation of scalar parameters for a system of evolutionary dissipative partial differential equations in which the state is partially observed. The algorithm takes advantage of the error that results from nudging a system with incorrect parameters with data from the true system. The intuitive nature of the algorithm makes its extension to several different systems immediate, and it allows for recovery of multiple parameters simultaneously.  We test the method on the Kuramoto-Sivashinsky equation in one dimension and demonstrate its efficacy in this context.
\end{abstract}



\section{Introduction}
Many mathematical models have parameters that must be identified correctly to match reality.  This is often a difficult task, as such parameters are difficult to measure and/or noise in the measurement process precludes identifying the parameters exactly.  Motivated by \cite{CaHuLa2020}, we develop a robust algorithm that estimates multiple scalar parameters in a partial differential equation (PDE) when those parameters can be written as coefficients of the various terms in the PDE.  As the current algorithm relies on the feedback control data assimilation mechanism introduced in \cite{azouani2014continuous}, we simultaneously identify not only the parameters of the system, but also obtain the dynamical state of the system based on reduced observations of the true state.

The method presented herein differs in important ways from other data-driven parameter recovery or equation discovery methods in the recent literature, the majority of which are \emph{offline} procedures in the sense that they seek to learn models from a previously collected set of training data.
The sparse identification of nonlinear dynamical systems (SINDy) framework \cite{Brunton2016sindy} generates a library of nonlinear candidate functions and leverages sparsity-promoting optimization techniques to generate the model that, as a linear combination of only a few of the candidate functions, best explains a set of observed dynamics. The methodology has also been extended to a parametric setting through periodic sparse updates to existing models \cite{MMKB2018sindyabrupt} or separating the parametric dependence from the state via group sparsity \cite{rudy2019data}.
Recent advances in deep learning propose neural network models that build the residual of a known governing equation into the training loss function, informing the optimization problem of the desired physics \cite{RPK2019pinns}.  Other deep learning approaches aim to directly approximate parameter-to-solution operators with neural networks \cite{LKALBS2020neuralOperator,LPK2021deeponet}, establishing a mapping between system parameters and solutions of the underlying PDEs. Such approaches massively overparameterize the system dynamics, endowing the resulting network with expressivity but often at the expense of training efficiency.
Bayesian statistics yields yet another framework wherein a posterior distribution with associated uncertainties can be constructed for the parameters of the PDE, or even the solution itself (see \cite{BoGlKr2020} for example).  The use of Markov chain Monte Carlo (MCMC) methods in this context provide robust methods for estimating the posterior distribution, and the advantages of the Bayesian approach are multifaceted \cite{DaSt2016}; however, such analyses are often computationally intensive as even the most advanced MCMC methods require thousands of evaluations of a numerical approximation to the PDE.

In contrast, the method we propose is a \emph{concurrent} (on-the-fly) procedure, providing updates to the unknown parameters of a target model simultaneous to the assimilation of measurements. We consider the setting in which the form of the governing equations is known, and therefore have no need to construct a library of candidate terms or overparametrize the solution representation, and later extend the approach to a setting where less is known about the governing equation. We leverage well-established methods originally developed to determine the true state of a system when the model is fully known, but in our setting we additionally infer the unknown system parameters.

Our derivation is motivated by the Azouani-Olson-Titi (AOT) \cite{azouani2014continuous} algorithm developed for the continuous assimilation of limited data into a known dissipative dynamical system.  Since \cite{azouani2014continuous}, this approach has been adapted to numerous other systems (see \cite{BeOlTi2015,BiHuLaPe2018,JoMaTi2017} for example) and settings in which the observed variables are various subsets of the full system's variables (e.g., \cite{FaJoJoTi2017,FaJoTi2015,farhat2016abridged,farhat2016charney}) or where different types of observations and nudging are utilized (see \cite{BiBrJo2020,LaPe2018} for example).  All of these extensions have broadened the applicability of the AOT algorithm but rely on the accurate representation of the underlying model parameters.  Recently, data assimilation for an imperfect model has been considered in different contexts \cite{CaHuLa2020,LeClMaBi2018,FaGlMaMcWh2020}.  It is rigorously established for two distinct settings that when a parameter of the system is unknown, the error in the assimilation will still converge up to the error in the unknown parameter (see \cite{CaHuLa2020,FaGlMaMcWh2020}). An upcoming article \cite{CaHuLaMaNgWh2021} takes these concepts a step further, deriving analogous update formulae for the three parameters of the Lorenz equations \cite{Lo1963} and carrying out a thorough numerical investigation to explore the limitations of these formulae.  In addition, asymptotic convergence to the true parameter values under suitable conditions is rigorously guaranteed, giving the first instance where a rigorous mathematical framework is available for this approach to learning parameters. We expand on all of this previous work, developing an algorithm that works not only for a single unknown parameter, but for multiple unknowns, and which theoretically applies to any dissipative dynamical system.

Our algorithm is derived in a general setting and is applicable to a wide class of dissipative PDEs. In this work, we examine the efficacy of the algorithm numerically for the Kuramoto-Sivashinsky equation (KSE).  KSE is selected because it presents chaotic spatial-temporal dynamics at a relatively low computational cost (see \cite{budanur2017unstable,cvitanovic2010state,edson2019lyapunov,wittenberg1999scale} for example studies of the complicated dynamics that arise in KSE).  In addition, KSE naturally provides the context for introducing several artificial parameters whose estimation is of significant import.

KSE has been derived in a variety of contexts, typically as a model for physical systems that are far from equilibrium.  The original derivation of KSE was to study instabilities in a reaction diffusion system  \cite{KuTs1975} and simultaneously on the instability and turbulent behavior of a single flame \cite{Si1977}.  KSE also appears in the study of plasmas \cite{CoKrTaRo1976,LaMaRuTa1975}, in irregular flow down an inclined plane \cite{SiMi1980}, and it is shown in \cite{MiVa1994} that KSE provides a general description for part of the dynamical evolution of systems following a certain type of bifurcation point.  Mathematically, KSE is often used as an ideal setting to investigate nonlinear dynamics, coherent structures, long-time behavior of solutions, and the testing of methods developed to investigate these properties (see \cite{CoEcEpSt1993,FoSeTi1989,GoFa2019,Gr2000,HyNi1986,NiScTe1985} for example).  Such earlier works have demonstrated that although KSE is simple and computationally cheap to simulate, it displays a wide variety of interesting dynamics similar to more complicated systems and is a reasonable testbed for the parameter estimation algorithm proposed below.
The paper is organized as follows. The derivation of the parameter learning algorithm is presented in
\Cref{sec:algorithm}, and \Cref{sec:numerics} presents numerical simulations applied to KSE that demonstrate the apparent numerical properties of the algorithm.
\Cref{sec:conclusions} draws conclusions and discusses extensions of these results to other systems and the potential for rigorous justification of the main algorithm.

\section{Derivation of the Algorithm}
\label{sec:algorithm}
In this section, we develop (but do not prove convergence of) a parameter recovery algorithm in a general setting.
Let $\Omega \subset \mathbb{R}^{d}$, $d\in\{1,2,3\}$, be an open set with Lipschitz continuous boundary $\partial \Omega$, and let $\mathcal{V} = L^{2}(\Omega)$ with symmetric inner product $\langle\cdot,\cdot\rangle$. All references to a vector space or properties of a vector space (orthogonality etc.) referred to below will apply to $L^2(\Omega)$ unless specified otherwise.

\subsection{Problem Statement}

Consider a PDE of the form
\begin{subequations}
\begin{align}
    \label{eq:u_evolution}
    \mathbf{u}_{t}
     + F(\mathbf{u})
     + \sum_{k=1}^n\lambda_{k} G_{k}(\mathbf{u})
    &= \mathbf{f},
\end{align}
where $\mathbf{u} = \mathbf{u}(\mathbf{x},t)$ is an unknown physical quantity (e.g., velocity), $\mathbf{u}_{t} = \frac{\partial}{\partial t}\mathbf{u}(\mathbf{x},t)$ is its partial derivative in time, $\mathbf{f} = \mathbf{f}(\mathbf{x},t)$ is a known forcing function,
$F$ and $G_{1},\ldots,G_{n}$ are (potentially nonlinear) spatial differential operators, and $\lambda_1, \ldots, \lambda_n\in\mathbb{R}$ are fixed real numbers.  Throughout the derivation we assume that solutions of \cref{eq:u_evolution} are classical \cite{Ev1998}, i.e., we will not concern ourselves with establishing the regularity of the solutions.

Our goal is to use partial observations of the state $\mathbf{u}$ to recover not only the full state, but also the unknown parameters $\lambda_{1},\ldots,\lambda_{n}$.  Before proceeding, we note that such a problem is not always well posed. For example, consider the following simplified, linear version of \eqref{eq:u_evolution},
\begin{align*}
    \mathbf{u}_t - \lambda_1\mathbf{u} +\lambda_2 A\mathbf{u}
    = \mathbf{f},
\end{align*}
where $A$ is a linear differential operator.  If the true solution $\mathbf{u}_\delta$ is a time-independent eigenfunction of $A$ with eigenvalue $\delta$, then the inverse problem reduces to attempting to solve
\begin{align*}
    (\lambda_2 \delta - \lambda_1)\mathbf{u}_\delta
    = \mathbf{f}
\end{align*}
for $\lambda_1$ and $\lambda_2$.
This is a degenerate situation wherein $\lambda_1$ and $\lambda_2$ cannot be uniquely identified by any algorithm. This issue can occur in other settings (e.g., with time-dependent solutions and/or nonlinear terms); however, the problem that we examine in \Cref{sec:numerics} does not suffer from lack of identifiability.

Motivated by \cite{CaHuLa2020,FaGlMaMcWh2020} and in tandem with \cref{eq:u_evolution}, we consider the coupled system
\begin{align}
    \label{eq:v_evolution}
    \mathbf{v}_{t} + F(\mathbf{v}) + \sum_{k=1}^n\widehat{\lambda}_{k}G_k(\mathbf{v})
    = \mathbf{f} + \mu (I_h(\mathbf{u})-I_h(\mathbf{v})),
\end{align}
\end{subequations}
where $\widehat{\lambda}_{k} = \widehat{\lambda}_{k}(t)$ is the estimate for $\lambda_{k}$ (which may evolve in time), $I_h$ is a linear interpolation operator corresponding to partial observations of the true system state, and $\mu > 0$ is called the nudging parameter.
The subscript $h$ refers to the level of resolution of the observation operator, i.e. $h\to 0$ would indicate observation of all relevant scales.
This is exactly the AOT formulation of \cite{azouani2014continuous}, except that $\widehat{\lambda}_k(t) \neq \lambda_k$, and these parameter estimates are updated in time along with the state estimate for $\mathbf{v}$. Our goal is to have $\mathbf{v}(t) \to \mathbf{u}(t)$ and $\widehat{\lambda}_{k}(t) \to \lambda_{k}$ up to an acceptable level of precision as $t$ increases.
Following arguments similar to those put forward in \cite{CaHuLa2020,FaGlMaMcWh2020}, under certain practical assumptions on \cref{eq:u_evolution} the convergence of $\mathbf{v}(t) \to \mathbf{u}(t) + \mathcal{O}(|\lambda_k-\widehat{\lambda}_k(t)|)$ is guaranteed; however, obtaining the desired parameter convergence $\widehat{\lambda}_{k}(t) \to \lambda_{k}$ requires an update formula for the parameters. This is the main subject of this work.

In \cite{CaHuLa2020}, two different update formulae are proposed for the 2D Navier-Stokes equations when the single unknown parameter $\lambda$ is the viscosity of the system.  Both approaches assume that $\mathbf{w} := \mathbf{u}-\mathbf{v}$ is sufficiently small so that quadratic terms in the error can be neglected, leading to an update formula for $\widehat{\lambda}$.  Such an assumption is valid, assuming that $\widehat{\lambda}$ is initially chosen close enough to $\lambda$ so that $\mathbf{w}$ is indeed small.  Critically, this update formula is applicable only when $n=1$, i.e., there is a single parameter to update.  The approach we take here does not suffer from either of these restrictions: we make no assumption on the magnitude of $\mathbf{w}$, and we can simultaneously recover multiple parameters at once.

\subsection{Parameter Estimation to Guide Synchronization}

We assume that we can only observe $I_h(\mathbf{u})$, rather than the full state $\mathbf{u}$. Furthermore, the time derivative $I_h(\mathbf{u}_{t})$ is not directly available (i.e., we do not have direct observations of $\mathbf{u}_{t}$) and must instead be numerically estimated via, e.g., finite differences of $I_h(\mathbf{u})$, which is feasible if we assume that $I_h$ is a linear operator (this assumption is carried throughout the following derivation).
Therefore, the metric of convergence that we use is the observable error, defined as $I_h(\mathbf{w}) = I_{h}(\mathbf{u}) - I_{h}(\mathbf{v})$.
While convergence of the observable error to zero is not sufficient to imply convergence of the parameters or full state estimates, it is necessary.
In other words, $\widehat{\lambda}_{k}(t)\to\lambda_{k}$ and $\mathbf{w}(t)\to\mathbf{0}$ can only occur simultaneously when $I_{h}(\mathbf{w}(t)) \to \mathbf{0}$.
To motivate the algorithm described here, we therefore examine the dynamics of $I_{h}(\mathbf{w})$.

\begin{remark}
We do not rigorously show that $I_h(\mathbf{w}(t))\to 0$ implies that $\mathbf{w}(t) \to 0$. However, traditional Foias-Prodi type estimates typically used to guarantee convergence of AOT-type algorithms apply in this setting with a correction term due to the differences in the true parameters and the estimated parameters.
Correction terms are considered in a different context in both \cite{CaHuLa2020} and \cite{FaGlMaMcWh2020}, although the current setting is more general.  The primary difficulty is that while such a correction term can certainly be estimated, it is updating in time, i.e., convergence of $\mathbf{w}(t) \to 0$ is dependent on the convergence of $|\lambda_k - \widehat{\lambda}_k(t)|\to 0$, and vice versa. That is, the dual convergence of the state and the parameters is not easily decoupled. This makes obtaining the usual error estimates challenging.
\end{remark}

Combining \eqref{eq:u_evolution} and \eqref{eq:v_evolution}, the evolution of the observable error can be written as
\begin{align*}
    \frac{d}{dt}I_{h}(\mathbf{w})
    &= I_{h}(\mathbf{u}_{t}) - I_{h}(\mathbf{v}_{t})
    \\
    &= I_h(\mathbf{u}_{t}) + I_{h}(F(\mathbf{v})) - I_h(\mathbf{f}) - \mu I_h(I_h(\mathbf{w})) + \sum_{k=1}^{n}\widehat{\lambda}_{k}I_{h}(G_{k}(\mathbf{v})),
\end{align*}
wherein we have used the linearity of $I_h$.
Defining
\begin{align*}
    \boldsymbol{\eta}
    = I_h(\mathbf{u}_{t}) &+ I_{h}(F(\mathbf{v})) - I_h(\mathbf{f})
    - \mu \left[ I_h(I_h(\mathbf{w})) - I_h(\mathbf{w})\right] + \sum_{k=1}^{n}\widehat{\lambda}_{k}I_{h}(G_{k}(\mathbf{v})),
\end{align*}
we have
\begin{align}
    \label{eq:ode-with-eta}
    \frac{d}{dt}I_{h}(\mathbf{w})
    = - \mu I_h(\mathbf{w}) + \boldsymbol{\eta}.
\end{align}
This leads to the following convergence result for $I_{h}(\mathbf{w})$.

\begin{proposition}
\label{thm:error-decay}
Assume that $\langle I_h(\mathbf{w}), \boldsymbol{\eta}\rangle = 0$ (the inner product is in $L^2$). Then $I_h(\mathbf{w}) \to \mathbf{0}$ exponentially in time with decay rate $\mu$.
\begin{proof}
Because $\langle\cdot,\cdot\rangle$ is symmetric and applying~\cref{eq:ode-with-eta},
\begin{align}
    \begin{aligned}
    \frac{1}{2}\frac{d}{dt}\left\|I_{h}(\mathbf{w})\right\|^2
    = \left\langle \frac{d}{dt}I_{h}(\mathbf{w}), I_{h}(\mathbf{w}) \right\rangle
&= -\mu \left\|I_h(\mathbf{w})\right\|^2 + \left\langle I_h(\mathbf{w}), \boldsymbol{\eta} \right\rangle,
    \end{aligned}
\end{align}
where
$
    \|\mathbf{u}\| = \sqrt{\langle\mathbf{u},\mathbf{u}\rangle}
$
is the norm induced by the inner product.
By the assumption that $\langle I_h(\mathbf{w}), \boldsymbol{\eta}\rangle = 0$, the previous equation simplifies to
\begin{align*}
    \frac{1}{2}\frac{d}{dt}\left\|I_h(\mathbf{w})\right\|^{2}
    &= -\mu \left\|I_h(\mathbf{w})\right\|^{2},
\end{align*}
which has solution $\left\|I_h(\mathbf{w}(t))\right\| = \left\|I_h(\mathbf{w}(0))\right\|e^{-\mu t}$.
\end{proof}
\end{proposition}

Note that $\boldsymbol{\eta}$ depends on time and the parameter estimates $\widehat{\lambda}_1,\ldots,\widehat{\lambda}_{n}$.
Our goal, then, is to choose $\widehat{\lambda}_1,\ldots,\widehat{\lambda}_{n}$ such that \Cref{thm:error-decay} applies and $I_{h}(\mathbf{w})\to\mathbf{0}$.
The requisite orthogonality condition is
\begin{align}
    \label{eq:scalar_constraint}
    \begin{aligned}
    0
    &= \langle I_h(\mathbf{w}), \boldsymbol{\eta}\rangle \phantom{\sum_{k=1}^{n}}
    \\
    &= \langle I_h(\mathbf{w}), I_h(\mathbf{u}_{t}) + I_{h}(F(\mathbf{v})) - I_h(\mathbf{f}) \rangle
    \\ &\qquad - \mu \langle I_h(\mathbf{w}), \left[ I_h(I_h(\mathbf{w})) - I_h(\mathbf{w})\right]\rangle
+ \sum_{k=1}^{n}\widehat{\lambda}_{k}\langle I_h(\mathbf{w}),  I_{h}(G_{k}(\mathbf{v}))\rangle.
    \end{aligned}
\end{align}
As long as at least one of the inner products $\langle I_h(\mathbf{w}), G_{k}(\mathbf{v}) \rangle$ is nonzero, we can choose $\widehat{\lambda}_{1},\ldots,\widehat{\lambda}_{n}$ so that \eqref{eq:scalar_constraint} is satisfied. This is an important assumption: if $I_h(\mathbf{w})$ happens to be orthogonal to $I_h(G_{k}(\mathbf{v}))$ for every $k=1,\ldots,n$, then the values of $\widehat{\lambda}_{1},\ldots,\widehat{\lambda}_{n}$ have no effect on $\left\langle I_{h}(\mathbf{w}), \boldsymbol{\eta}\right\rangle$.

In the case of a single parameter ($n = 1$), \eqref{eq:scalar_constraint} uniquely determines the parameter estimate $\widehat{\lambda}_{1}$, provided $\langle I_h(\mathbf{w}), G_{1}(\mathbf{v}) \rangle \neq 0$. However, if there are multiple parameters ($n > 1$), \eqref{eq:scalar_constraint} is underdetermined, which prompts us to impose additional constraints. Fix $\mathbf{e}_1$ = $I_h(\mathbf{w})$ as in \cref{eq:scalar_constraint} and let $\mathcal{G}(\mathbf{v}) = \textrm{span}(\{I_h(G_k(\mathbf{v}))\}_{k=1}^{n})$.
If $\mathcal{G}(\mathbf{v})$ has dimension $n$ (in other words, if the set $\{I_h(G_k(\mathbf{v}))\}_{k=1}^{n}$ is linearly independent), then we can select $n-1$ additional linearly independent functions $\mathbf{e}_2, \ldots, \mathbf{e}_n \in \mathcal{G}(\mathbf{v})$ that are orthogonal to $\boldsymbol{\eta}$, i.e.,
$\langle \mathbf{e}_i, \boldsymbol{\eta}\rangle = 0$ for $i=1,\ldots,n$.
This requirement gives rise to the $n\times n$ system of equations
\begin{align}
    \nonumber &&
    A_{i,k} &= \left\langle \mathbf{e}_{i}, G_{k}(I_{h}(\mathbf{v}))\right\rangle,
    \\ \label{eq:proxy_multi_param_sys}
    \mathbf{A} \widehat{\boldsymbol{\lambda}}
    &= \mathbf{b},
    &
    \widehat{\boldsymbol{\lambda}}
    &= [~\widehat{\lambda}_{1},~\ldots,~\widehat{\lambda}_{n}~]^{\mathsf{T}},
    \\ \nonumber &&
    b_{i}
    &= \left\langle \mathbf{e}_{i}, I_{h}(\mathbf{f} - \mathbf{u}_{t} - \mathbf{F}(\mathbf{v})) + \mu[I_{h}(I_{h}(\mathbf{w})) - I_{h}(\mathbf{w})] \right\rangle,
\end{align}
wherein we have suppressed the dependence on time for simplicity.
The assumption that $\text{dim}(\mathcal{G}(\mathbf{v})) = n$ guarantees that this system will have a unique solution.

\begin{remark}
The particular choice of the basis functions $\mathbf{e}_{1},\ldots,\mathbf{e}_{n}$ introduces both theoretical and numerical/practical considerations.
\Cref{thm:error-decay} advocates setting $\mathbf{e}_1 = I_h(\mathbf{w})$ as in \cref{eq:scalar_constraint}. Theoretically, it may seem advantageous to choose $\mathbf{e}_2$ so that we can guarantee some higher order norm of $I_h(\mathbf{w})$ decays exponentially, much as selecting $\mathbf{e}_1 = I_h(\mathbf{w})$ leads to decay of $I_h(\mathbf{w})$ in $L^2(\Omega)$.  However, it is not clear how to do this in practice and can lead to a numerically expensive calculation, so we instead simply choose the $\mathbf{e}_{2},\ldots,\mathbf{e}_{n}$ to form an orthonormal set to obtain a numerically well-conditioned linear system.
\end{remark}

Suppose now that $\mathbf{w}$ is small, that is, $\mathbf{v} \approx \mathbf{u}$. In this case, we have
\begin{align*}
    \boldsymbol{\eta}
    \approx \tilde{\boldsymbol{\eta}}
    :=&~I_h(\mathbf{u}_{t}) + I_{h}(F(\mathbf{u})) - I_h(\mathbf{f}) + \sum_{k=1}^{n}\widehat{\lambda}_{k}I_{h}(G_{k}(\mathbf{u}))
    \\
    =&~I_h\Big(\mathbf{u}_{t} + F(\mathbf{u}) - \mathbf{f} + \sum_{k=1}^{n}\widehat{\lambda}_{k}G_{k}(\mathbf{u})\Big),
\end{align*}
provided the operators $F$ and $G_1,\ldots,G_{n}$ are sufficiently smooth.\footnote{This is a continuity assumption: if $\mathbf{v}\approx\mathbf{u}$, then $F(\mathbf{v})\approx F(\mathbf{u})$ and $G_{k}(\mathbf{v})\approx G_{k}(\mathbf{u})$, $k=1,\ldots,n$.}
Observe that if the parameter estimates $\widehat{\lambda}_{1},\ldots,\widehat{\lambda}_{n}$ are replaced with the true parameters $\lambda_1,\ldots,\lambda_{n}$, then $\tilde{\boldsymbol{\eta}} \equiv \mathbf{0}$ by \eqref{eq:u_evolution}. This suggests that parameter estimates obtained by solving \eqref{eq:proxy_multi_param_sys} are reasonable approximations to the true parameters when the estimated state $\mathbf{v}$ is close to the true state $\mathbf{u}$.

Thus, under reasonable assumptions, parameter estimates given by the solution of~\cref{eq:proxy_multi_param_sys} with appropriately selected basis functions $\mathbf{e}_{1},\ldots,\mathbf{e}_{n}$ both drive the observable error to zero and at least approximately recover the true parameters when the state estimate converges. These properties are not sufficient to guarantee convergence of the parameter estimates; nevertheless, we will demonstrate that a numerical estimation algorithm based on \cref{eq:proxy_multi_param_sys} is effective in practice.

\subsection{Numerical Considerations}

As described above, we select $\mathbf{e}_1$ = $I_h(\mathbf{w})$ and choose the other $\mathbf{e}_i$ to be orthonormal, ensuring that \cref{eq:proxy_multi_param_sys} is well conditioned. The straightforward approach is then to solve \cref{eq:proxy_multi_param_sys} at each time step and use the obtained parameter estimates to advance $\mathbf{v}$. However, this strategy can lead to rapid/discontinuous changes in the parameter at each time step.  Such discontinuous adjustments quickly produce numerical instabilities that are unavoidable if the update is applied at each time step, particularly when the initial parameter guess is not close to the true value.

To avoid discontinuous updates of the parameter estimates, we numerically integrate
\begin{equation}
    \label{eq:relaxation_ode}
    \frac{d \widehat{\lambda}_k(t)}{dt}
    = -\alpha (\widehat{\lambda}_k(t) - \widetilde{\lambda}_k(t)),
\end{equation}
where $\widetilde{\lambda}_k(t)$ is obtained by solving \cref{eq:proxy_multi_param_sys} at each time $t$ and $\alpha > 0$ is a relaxation parameter.  We then solve the coupled system given by \cref{eq:u_evolution}--\cref{eq:v_evolution} and \cref{eq:relaxation_ode}.  This gives a parameter update that avoids discontinuous changes in the parameters of the system, hence preserving the rigorously justified smoothness of the solutions $\mathbf{v}$ and $\mathbf{u}$.  The effect of the relaxation are examined in \cref{fig:relaxation} for KSE as described in \cref{sec:smoothing}.

\Cref{alg:parameter-recovery} details the entire parameter recovery procedure.
Note that the algorithm is \emph{concurrent} in the sense that the state estimates $\mathbf{v}_{j} := \mathbf{v}(t_j)$ can be computed as soon as the projected data at the corresponding time, $\widehat{\mathbf{u}}_{j} := I_{h}(\mathbf{u}(t_{j}))$, becomes available.
We denote this in the algorithm with the observation operator $\mathcal{Q}: t\mapsto I_{h}(\mathbf{u}(t))$.
The state estimates are updated by time stepping \cref{eq:v_evolution}, which we denote with a solution operator $\mathcal{S}: (\mathbf{v}_{j},\widehat{\lambda}_{1},\ldots,\widehat{\lambda}_{n}) \mapsto \mathbf{v}_{j+1}$ with uniform time step $\delta t = t_{j+1} - t_{j}$.
We make the simplifying assumption that $I_{h}$ is idempotent, i.e., $I_{h}(I_{h}(\mathbf{w})) = I_{h}(\mathbf{w})$, which slightly simplifies step~\ref{step:b-entries} of the algorithm compared with \cref{eq:scalar_constraint}--\cref{eq:proxy_multi_param_sys}.
The operators $I_{h}$ used in our numerical experiments satisfy this property.

Before presenting results of this parameter estimation strategy on KSE, we emphasize that the approach depends critically on assimilating data continuously in time so that the time derivative of the observed state $I_h(\mathbf{u}_t)$ can be approximated accurately.  As shown later in \Cref{fig:convergence_order}, the error in the parameter approximation is directly dictated by this finite difference approximation, so that discrete measurements of $I_h(\mathbf{u})$ that are spaced too far apart in time will cause the parameter estimation to break down.  All numerical results reported in \Cref{sec:numerics} are based on continuous data assimilation in the sense that $I_h(\mathbf{u})$ is identified at every time step of the numerical solver.

\begin{algorithm}[ht]
\begin{algorithmic}[1]
\Procedure{MPE}{
    Initial time $t_0\in\mathbb{R}$,
    time step $\delta t > 0$,
    projected state observer $\mathcal{Q}:\mathbb{R}\to\mathcal{V}$,
    PDE solver $\mathcal{S}:\mathcal{V}\times\mathbb{R}^{n}\to\mathcal{V}$ (time stepper),
    initial parameter guesses $\widehat{\lambda}_{1},\ldots,\widehat{\lambda}_{n}\in\mathbb{R}$,
    initial state estimate $\mathbf{v}_0 \in \mathcal{V}$,
    finite difference order $p\in\mathbb{N}$,
    projector $I_{h}:\mathcal{V}\to\mathcal{V}$,
    nudging parameter $\mu > 0$,
    relaxation parameter $\alpha > 0$
    }
    \vspace{.0625in}
    \State $t \gets t_0$
    \For{$j=0,1,\ldots$}
        \LineComment{Observe projected state and estimate its time derivative.}
        \State $\widehat{\mathbf{u}}_{j} \gets \mathcal{Q}(t)$
\Comment{Projected state observation.}
\State $\widehat{\mathbf{u}}'_{j} \gets \textproc{BDF}(\widehat{\mathbf{u}}_{j}, \ldots, \widehat{\mathbf{u}}_{j-p})$
            \Comment{Backwards difference formula.}

        \LineComment{Construct and orthonormalize the basis vectors.}
        \State $\mathbf{e}_1 \gets \widehat{\mathbf{u}}_{j} - I_{h}(\mathbf{v}_{j})$
        \For{$i=2, \ldots,n$}
        	\State $\mathbf{e}_{i} \gets I_{h}(G_{i-1}(\mathbf{v}_{j}))$
        \EndFor
        \State $\mathbf{e}_1, \ldots, \mathbf{e}_n \gets \textproc{MGS}(\mathbf{e}_1, \ldots, \mathbf{e}_n)$
            \Comment{Modified Gram-Schmidt.}

        \LineComment{Construct system for estimating parameters.}
        \State $\mathbf{b} \in \mathbb{R}^{n}$, $\mathbf{A}\in\mathbb{R}^{n\times n}$
            \Comment{Initialize empty n x n linear system.}
        \For{$i = 1, \ldots, n$}
            \State $b_{i} \gets \langle \mathbf{e}_{i}, I_{h}(\mathbf{f}) - \widehat{\mathbf{u}}_{j}' - I_{h}(F(\mathbf{v}_{j}))\rangle$
            \label{step:b-entries}
            \For{$k=1, \ldots, n$}
 \State $A_{i,k} \gets \langle \mathbf{e}_{i}, I_{h}(G_{k}(\mathbf{v}_{j}))\rangle$
            \EndFor
        \EndFor

		\LineComment{Get point-in-time estimates and smooth.}
		\State $\widetilde{\lambda}_{1},\cdots,\widetilde{\lambda}_{n} \gets \mathbf{A}^{-1}\mathbf{b}$
        \For{$k=1, \ldots, n$}
            \State $\widehat{\lambda}_{k} \gets \left(1 + \frac{1}{2}\alpha \delta t\right)^{-1}\left[\left(1 - \frac{1}{2}\alpha \delta t\right) \widehat{\lambda}_{k} + \alpha\delta t \widetilde{\lambda}_{k}\right]$
        \EndFor

        \LineComment{Update estimated state with AOT nudging.}
        \State $\mathbf{v}_{j+1} \gets \mathcal{S}(\mathbf{v}_{j}, \widehat{\lambda}_1, \ldots, \widehat{\lambda}_n) + \mu \delta t (\widehat{\mathbf{u}}_{j} - I_{h}(\mathbf{v}_{j}))$
        \State $t = t + \delta t$
\label{line:time-step-true-state}
	\EndFor

    \State \textbf{return} $\widehat{\lambda}_{1},\ldots,\widehat{\lambda}_{n}$
\EndProcedure

\end{algorithmic}
\caption{Concurrent Multi-Parameter Estimation}
\label{alg:parameter-recovery}
\end{algorithm}

\section{Numerical Results}
\label{sec:numerics}
The non-dimensionalization of KSE that we consider has the following form:
\begin{equation}
    \label{eq:kse_def}
    u_t + u_{xxxx} + \lambda u_{xx} + u u_x = 0.
\end{equation}
We impose periodic boundary conditions on the spatial domain $[0, 2\pi L)$. The chaos of the system depends on both $L$ and $\lambda$. Higher values of $L$ or $\lambda$ increase the exhibited chaos \cite{cvitanovic2010state,edson2019lyapunov,HyNi1986,wittenberg1999scale}. We choose $L$ = 16, $\lambda$ = 1 for all simulations as this represents a sufficiently chaotic state \cite{KaTr2005}.
In all simulations, we use a Fourier collocation scheme with a fourth-order Runge-Kutta implicit-explicit time-stepper \cite{KENNEDY2003139}.  The AOT nudging is applied via a forward Euler time-stepping scheme to alleviate some concerns over the implicit nature of the RK method,\footnote{As noted in \cite{OlTi2008} and observed further in \cite{FaGlMaMcWh2020}, multi-stage methods such as RK effectively interpolate the data between time steps, introducing some minor interpolation error in the nudging and hence the parameter estimation as well.} and the parameter smoothing \cref{eq:relaxation_ode} is integrated via the trapezoid method. The number of Fourier collocation points is a minimum of 512 in each simulation, which guarantees that even the smallest scales are sufficiently resolved. In most of our experiments, the observation projection operator is a Galerkin truncation onto the lowest 21 modes, meaning that $h \approx \frac{32\pi}{21}$ for all computations reported below.  All simulations are run out from $t=0$ to at least $t_f=50$ in non-dimensional units.  The code for running the simulations and generating the figures in this paper is publicly available at \texttt{\url{https://github.com/Parametric-Data-Assimilation/kse-multiparameter-learning}}.

\subsection{Single-parameter Learning}

We begin with a discussion of the performance of \cref{alg:parameter-recovery} when learning a single unknown parameter $\lambda$, i.e.,
\begin{subequations}
\begin{align}
    \label{eq:kse-singleparam-truth}
    u_t + u_{xxxx} + \lambda u_{xx} + u u_x
    &= 0, &&\text{(truth)}
    \\
    \label{eq:kse-singleparam-assim}
    v_t + v_{xxxx} + \widehat{\lambda} v_{xx} + v v_x
    &= \mu(I_{h}(u) - I_{h}(v)),
    &&\text{(assim)}
\end{align}
\end{subequations}
denoting $v$ for the nudging variable.

\subsubsection{Effect of the Relaxation and Nudging Parameters}
\label{sec:smoothing}

\begin{figure}
    \centering
    \includegraphics[width=\textwidth]{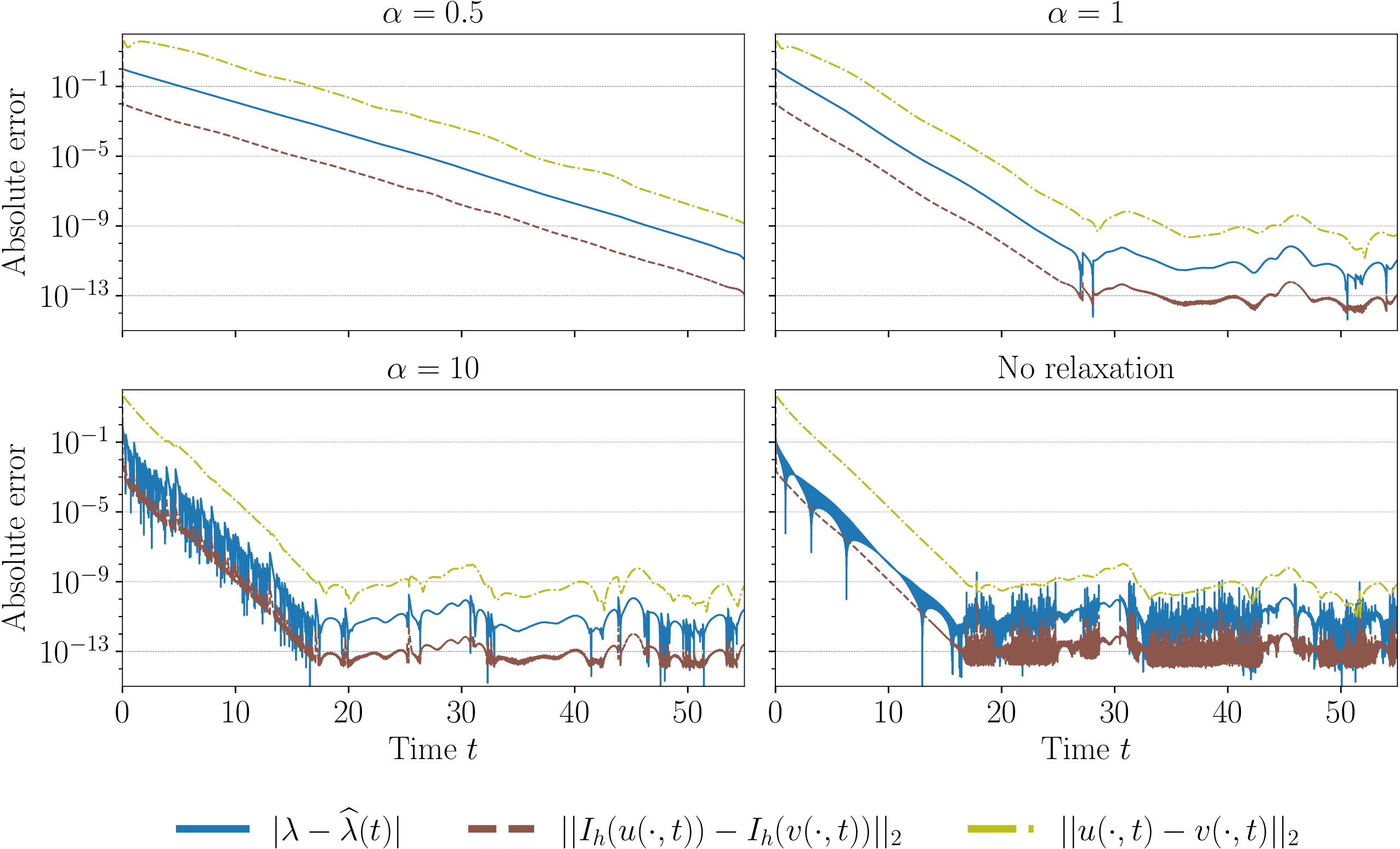}
    \vspace{-0.5cm}
    \caption{Finite relaxation as in \cref{eq:relaxation_ode} affects the rate and quality of convergence: strong relaxation (top row) smooths the convergence significantly but slightly decreases the convergence rate (see also \cref{fig:alpha_convergence_rate}); mild relaxation (bottom left) preserves the convergence rate that occurs when relaxation is absent (bottom right) but with some smoothing after convergence.}
    \label{fig:relaxation}
\end{figure}

We generally observe exponential decay of the projected error in time as demonstrated in \cref{fig:relaxation}. Defining $w := u - v$, both the solution error $\|w(\cdot,t)\|_2$ and the parameter estimate error $|\lambda-\widehat{\lambda}(t)|$ exhibit corresponding exponential decay in time, where $\|\cdot\|_{2}$ is the $L^2$ norm over the spatial domain.
The decay continues until the error is saturated by the effects of the finite difference time derivative approximation (discussed in \cref{sec:fd-order}), at which point the errors flatten out. Until then, we observe
\begin{align*}
    \|w(\cdot,t)\|_{2} \approx e^{-\beta t} \|w(\cdot,0)\|_{2},
\end{align*}
where $\beta > 0$ is the convergence rate. We expect $\beta$ to depend on the nudging parameter $\mu$ of \cref{eq:kse-singleparam-assim} and the relaxation parameter $\alpha$ of \cref{eq:relaxation_ode} used in the parameter update. In theory, increasing either $\mu$ or $\alpha$ should lead to faster convergence. However, extreme values of either parameter can violate the numerical stability criteria of the selected time-stepping scheme. The natural scaling for $\mu$ is with $\delta t^{-1}$ (the inverse of the time step). The stability region of a simple nudging as applied here indicates that $\mu \geq 2 \delta t^{-1}$ leads to instability, so we take $\mu$ = $1.8 \delta t^{-1}$ as the default nudging parameter.

\begin{figure}
    \centering
    \includegraphics[width=\textwidth]{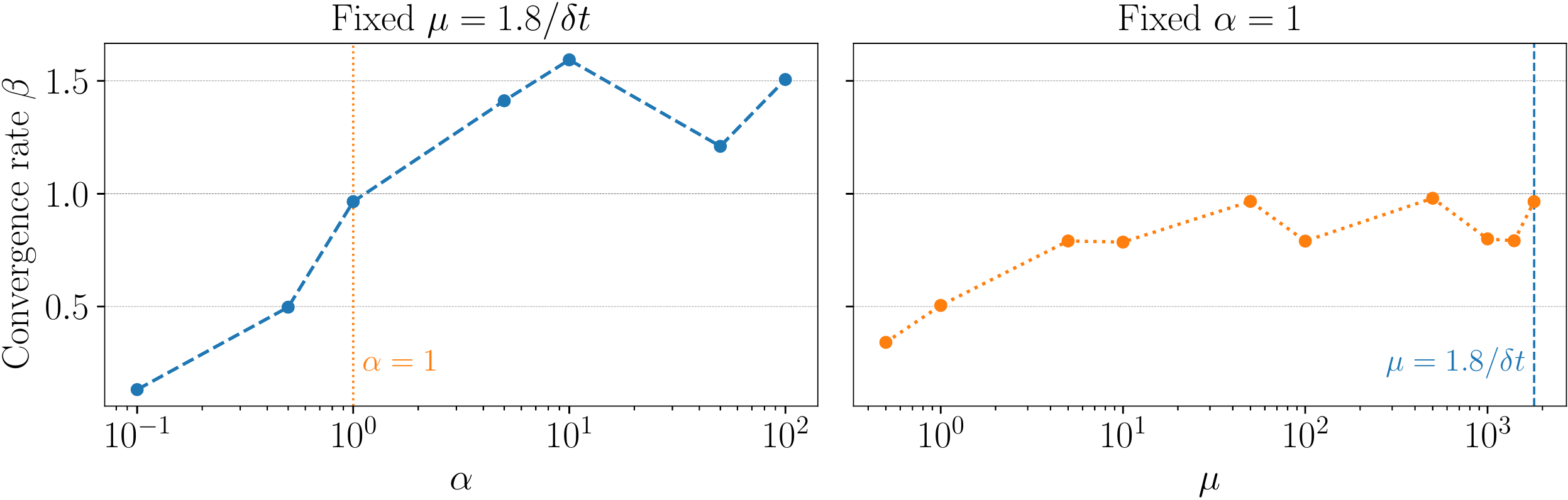}
    \vspace{-0.5cm}
    \caption{The convergence rate depends on the smoothing parameter $\alpha$ (left) and the assimilation constant $\mu$ (right). However, faster convergence is accompanied by non-smooth decay of the error in the parameter estimate (see \cref{fig:relaxation}).}
    \label{fig:alpha_convergence_rate}
\end{figure}

To explore the effect of $\alpha$ on convergence, we fix $\mu = 1.8 \delta t^{-1}$ and, for varying values of $\alpha$, estimate the convergence rate $\beta$ in two steps: 1) define the time of convergence $t_c$ as the first point in time when the error $\|w(\cdot,t_c)\|_2$ is no better than the error $\|w(\cdot,t_f)\|_2$ at the final time $t_f$, then 2) estimate $\beta$ from the time of convergence $t_c$ and the logarithmic ratio of the initial error and the error at $t_c$, i.e.,
\begin{align*}
    \beta \approx -\frac{\log\|w(\cdot,t_c)\|_2 - \log\|w(\cdot,t_0)\|_2}{t_{c} - t_{0}}.
\end{align*}
\Cref{fig:alpha_convergence_rate} shows that $\beta$ generally increases with $\alpha$ until about $\alpha = 10$. Repeating the experiment with fixed $\alpha = 1$ and varying $\mu$ shows that scaling $\mu$ down to $\mu < \mathcal{O}(10)$ decreases the convergence rate, but the dependence of $\beta$ on $\mu$ is not as significant as the dependence on $\alpha$.
In addition to affecting the convergence rate, $\alpha$ has a strong effect on the smoothness of the convergence, though this is more difficult to precisely quantify. \Cref{fig:relaxation} shows the convergence for a few choices of $\alpha$, with $\mu = 1.8\delta t^{-1}$. For single parameter estimation, spikes in the convergence without relaxation as seen in the bottom right plot of \cref{fig:relaxation} are not concerning, but such discontinuous jumps in the parameters can be problematic when estimating multiple parameters. Based on these observations, we set $\alpha = 1$ and $\mu = 1.8\delta t^{-1}$ for all subsequent experiments as a compromise between the speed and the smoothness of convergence.

\subsubsection{Choice of Observation Operator}
\label{sec:projection}

\begin{figure}
    \centering
    \includegraphics[width=\textwidth]{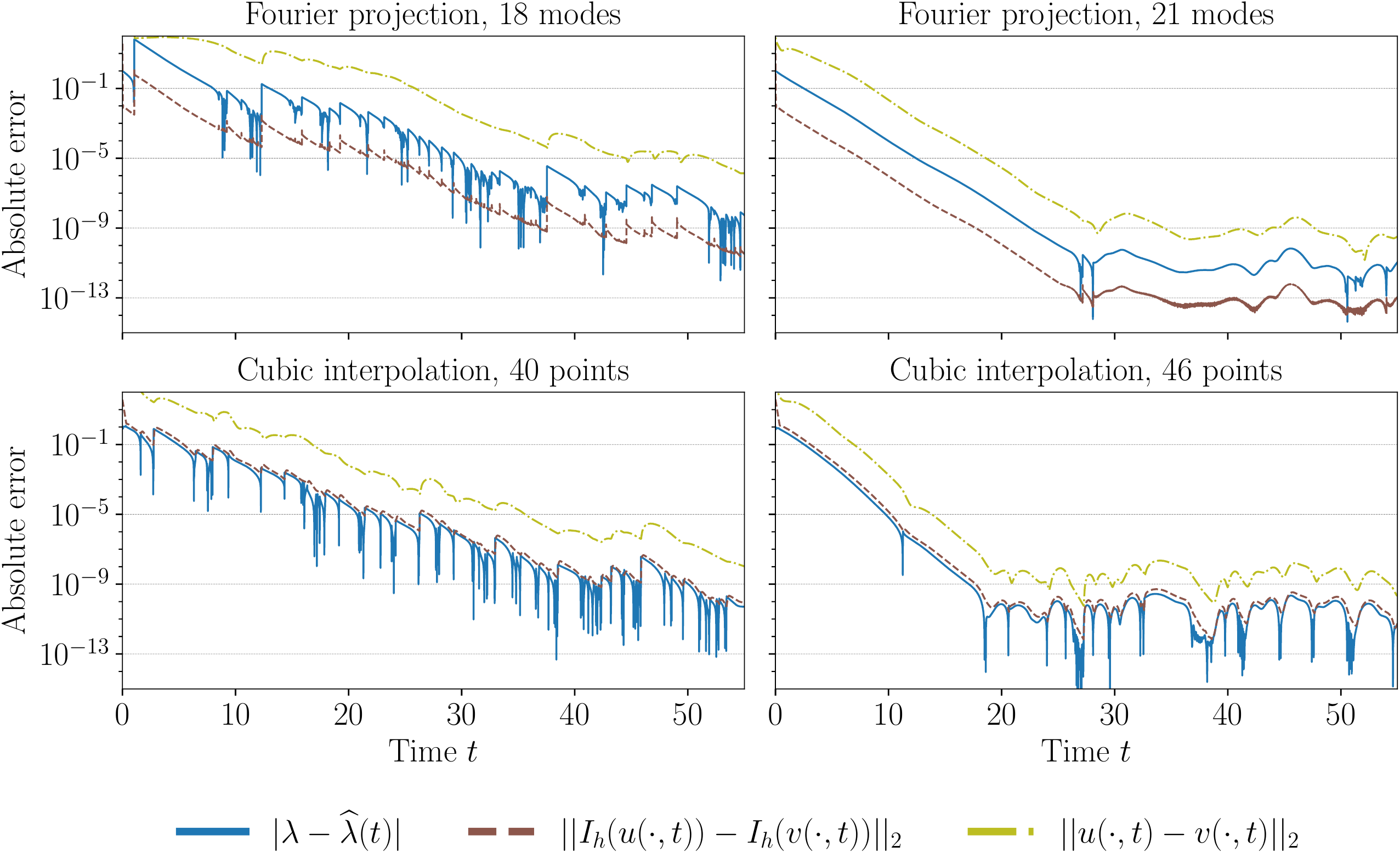}
    \vspace{-0.5cm}
    \caption{Convergence of single-parameter learning with Fourier truncation for the projection operator $I_h$ (top) and a cubic spline interpolator for $I_h$ (bottom). The plots on the left represent the minimal number of observations required for synchronization in each strategy.}
    \label{fig:interpolation}
\end{figure}

Up to this point, for the observation operator $I_h$ our experiments use a Fourier projection that truncates the state to the lowest $21$ Fourier modes.
This particular choice is natural given that we take a spectral approach to numerically solving \cref{eq:kse-singleparam-truth}--\cref{eq:kse-singleparam-assim} (see \cref{apdx:solver}).
\Cref{fig:interpolation} shows results comparable to \cref{fig:relaxation} but with only $18$ observed Fourier modes, which exhibits convergence but with less smoothness and at a slightly slower rate. Using fewer Fourier modes than 18 does not result in synchronization.

To further validate our method for a wider class of physical settings, we experiment with observation operators $I_h$ based on point-wise observations at evenly spaced points in the domain. Specifically, \cref{fig:interpolation} shows results when using cubic spline interpolation with $40$ and $46$ observed points. While not shown, using a linear or quadratic interpolation for $I_h$ yields similar results. The default nudging parameter $\mu$=$1.8\delta t ^{-1}$ does not result in synchronization in these experiments, and consequently we use the slightly more conservative value $\mu = 10$ in conjunction with point-wise observations. Even so, the convergence in the point-wise observations case is similar to the Fourier projection case, although the number of point-wise observations required to obtain convergence (40) is just over twice the number of Fourier modes needed (18).
At these marginal observation levels the convergence is rough, although adding just a few more observations results in a much smoother synchronization.  On the other hand, significantly increasing the number of observations beyond this marginal threshold does not substantially improve either the final total error, or the convergence rate.

\subsubsection{Order of Accuracy}
\label{sec:fd-order}

\begin{figure}
    \centering
    \includegraphics[width=\textwidth]{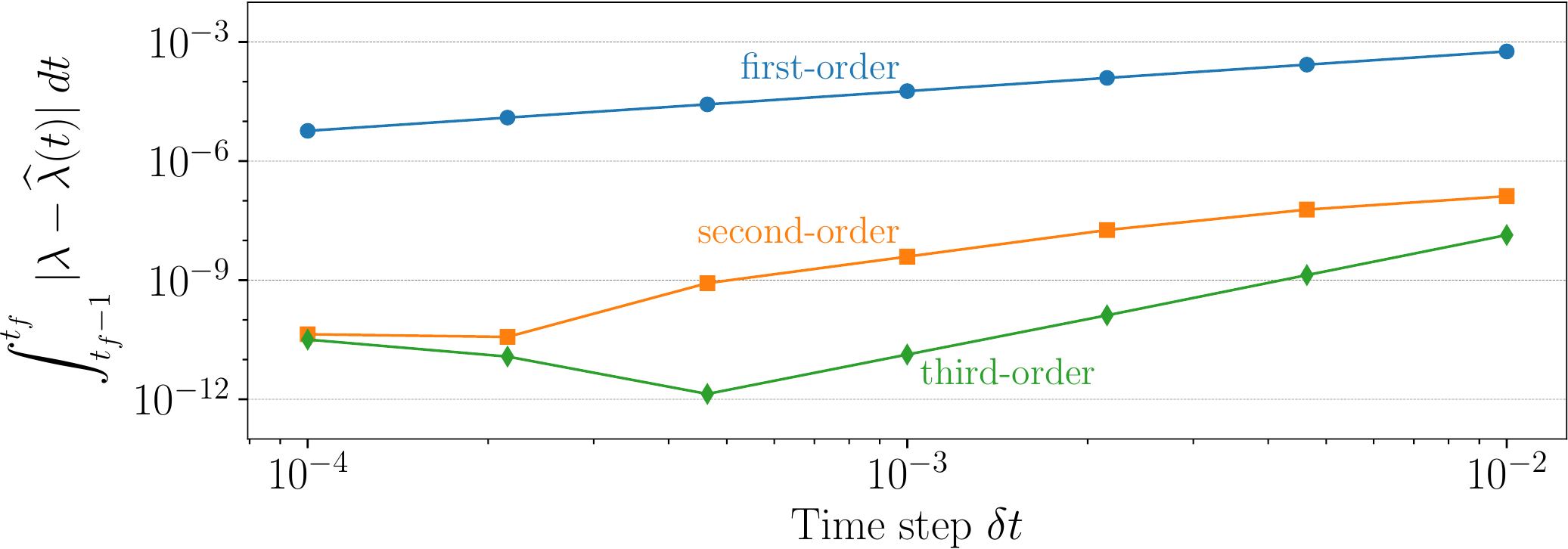}
    \vspace{-0.5cm}
    \caption{Error in the parameter estimate with respect to the time step $\delta t$ using first-order, second-order, and third-order backward finite difference schemes to estimate the projected time derivative $I_h(u_t)$. The order of the error is $1.0004$, $1.9489$, and $3.0023$, respectively (for the third-order method we exclude the two leftmost points). Excluding the leftmost point for the second-order finite differences raises the order estimate to $2.0535$.}
    \label{fig:convergence_order}
\end{figure}

The accuracy of the parameter estimate, and hence the estimate of the solution $u$, is limited by three factors: 1) the error introduced by the numerical scheme for the PDE, 2) round-off error, and 3) errors introduced by finite difference estimates of $I_h(u_t)$. Of these factors, the finite difference error for estimating $I_h(u_t)$ is the most significant, and therefore the order of the corresponding finite difference scheme controls the overall accuracy of the final parameter estimate with respect to the time step $\delta t$. \Cref{fig:convergence_order} shows the average error in the parameter estimate over the final simulation second
as a function of $\delta t$, with results using first-order, second-order, and third-order backward difference methods to estimate $I_h(u_t)$ (we use only backward differences as only data at the current and previous time steps are available).
We note that errors from the underlying PDE solver, and particularly the data assimilation aspect of the time stepping with the Runge-Kutta scheme,\footnote{See \cite{FaGlMaMcWh2020} for a description of this kind of error in a different setting, and \cite{OlTi2008} for a more thorough discussion of this issue.} limit the overall error to approximately $10^{-12}$, which is seen by the saturation of the error for both the second- and third-order finite difference methods in \cref{fig:convergence_order}. Omitting these final points, we obtain a linear regressed exponent that fits nearly exactly with the order of the finite difference method for the temporal derivative.  As these numerical investigations are meant to provide a proof of concept only, the error introduced via the RK time stepping is an acceptable compromise to facilitate ease of implementation and accessibility of the code.

We consider the propagation of the error in the third-order backward difference method specifically to understand the increase in the parameter error for $\delta t < 10^{-3}$.  If we allow $\xi$ to be the error from the Runge-Kutta time-stepping in the data assimilation step, then the numerically `converged' state of the system should be $\widetilde{I}_h(u) \approx I_h(u) + \xi$, where $\xi = \xi(\delta t,N,\epsilon_\textup{machine})$ is a function of the time step size, number of collocation points $N$, and machine epsilon $\epsilon_\textup{machine}$.  Using a third-order backward difference to approximate $I_h(u_t)$ in this case leads to
\begin{equation*}
    \widetilde{I}_h(u_t) \approx I_h(u_t) + \frac{\xi}{\delta t} + \mathcal{O}\left(\delta t^3\right).
\end{equation*}
The issue is that, as noted in \cite{FaGlMaMcWh2020} in a different setting, the multi-stage Runge-Kutta method introduces an error $\xi > \epsilon_\textup{machine}$ independent of $\delta t$ itself.  Thus, there is a certain point at which decreasing $\delta t$ counterintuitively leads to increased error in $\widetilde{I}_h(u_t)$, and hence an increased error in the parameter estimation.  The rub of the matter is that the two smallest values of $\delta t$ for the third-order difference in \cref{fig:convergence_order} are purely artifacts of the RK time stepping, and are hence numerically irrelevant to a discussion on the order of convergence.

\subsection{Multi-parameter Learning}

The primary contribution of \cref{alg:parameter-recovery} is the ability to learn multiple unknown parameters in a single equation in a concurrent fashion. To this end, we consider a \textit{generalized KSE} for the assimilating state $v(x,t)$, defined as
\begin{subequations}
\begin{align}
    \label{eq:kse_standard}
    u_t + \lambda_{4}u_{xxxx} + \lambda_{2}u_{xx} + \lambda_{5}u u_x
    &= 0, &&\text{(truth)}
    \\ \label{eq:kse_def_generalized}
    v_t + \widehat{\lambda}_{4}v_{xxxx} + \widehat{\lambda}_{3}v_{xxx} + \widehat{\lambda}_{2} v_{xx} + \widehat{\lambda}_{1}v_{x} + \widehat{\lambda}_{5}v v_x
    &= \mu(I_{h}(u) - I_{h}(v)).
    &&\text{(assim)}
\end{align}
\end{subequations}
Once again we consider $x\in[0,2\pi L)$ with $L=16$ fixed and run all simulations from time $t=0$ to $t_f=50$.  Following the discussion above for the single parameter case, we set $\delta t=10^{-3}$, $\mu = 1.8\delta t^{-1}$, $\alpha = 1$, define the projection operator $I_h$ as the Fourier projection that truncates to the lowest $21$ modes, and use third-order backward differences to approximate the projection of the temporal derivative $I_h(u_t)$.

The goal of this investigation is not only to identify the correct parameters of the specified PDE, but also to discover the structure of an unknown equation.  Hence, we have expanded the standard KSE to include unknown parameters on all terms and additional parameters on linear odd-order derivatives.
The true KSE corresponds to the case when $\lambda_1=\lambda_3=0$ and $\lambda_2=\lambda_4=\lambda_5=1$. In each experiment, we choose a subset of the $\lambda_k$ to learn, and set the rest to their true values. The simulation is initialized with an initial guess of $\widehat{\lambda}_k=2$ for every $\lambda_k$ that is unknown.

We focus here on identifying the coefficient $\lambda_5$ in front of the known quadratic nonlinearity $uu_x$, although the derivation of \cref{alg:parameter-recovery} presented above technically holds for additional nonlinear terms.  The reason for this restriction in the current setting is that adding other nonlinear terms to KSE destabilizes the dynamical evolution so that neither the original system nor the nudged system can be simulated easily.  Hence, the restriction to a single nonlinear term here is more a product of the specific dynamical system we are investigating (KSE) than a product of the proposed algorithm.

We preface our numerical results by first discussing a few implementation details which are not obvious from the theoretical derivation of the algorithm.

\subsubsection{Choice of Basis}

A crucial piece of the algorithm is the choice of basis functions $e_1,\ldots,e_N$. Other than setting $e_1 = I_h(w)$ (to ensure exponential decay of $I_h(w)$ in $L^2$), it is not clear how to choose these functions. For numerical efficiency, we would like to use information that has already been computed. A natural choice is then to use the $I_h(G_k(v))$, as these are computed in the simulation itself. Here another ambiguity arises---we need only $n-1$ more basis elements, but there are $n$ of the $I_h(G_k(v))$, hence we may exclude one. Most of the time, the choice does not affect convergence. However, we find that when $G_k$ corresponds to the nonlinear term in the generalized KSE, including it as a basis element causes the algorithm to not converge. Consequently, we always exclude this term from the basis.

For numerical stability, the basis vectors are orthonormalized in $L^2$ before constructing the system  \eqref{eq:proxy_multi_param_sys}.  While it is technically possible for the system to be degenerate, we do not observe this in practice at least in this setting.

\subsubsection{Convergence}

\begin{figure}
    \centering
    \includegraphics[width=\textwidth]{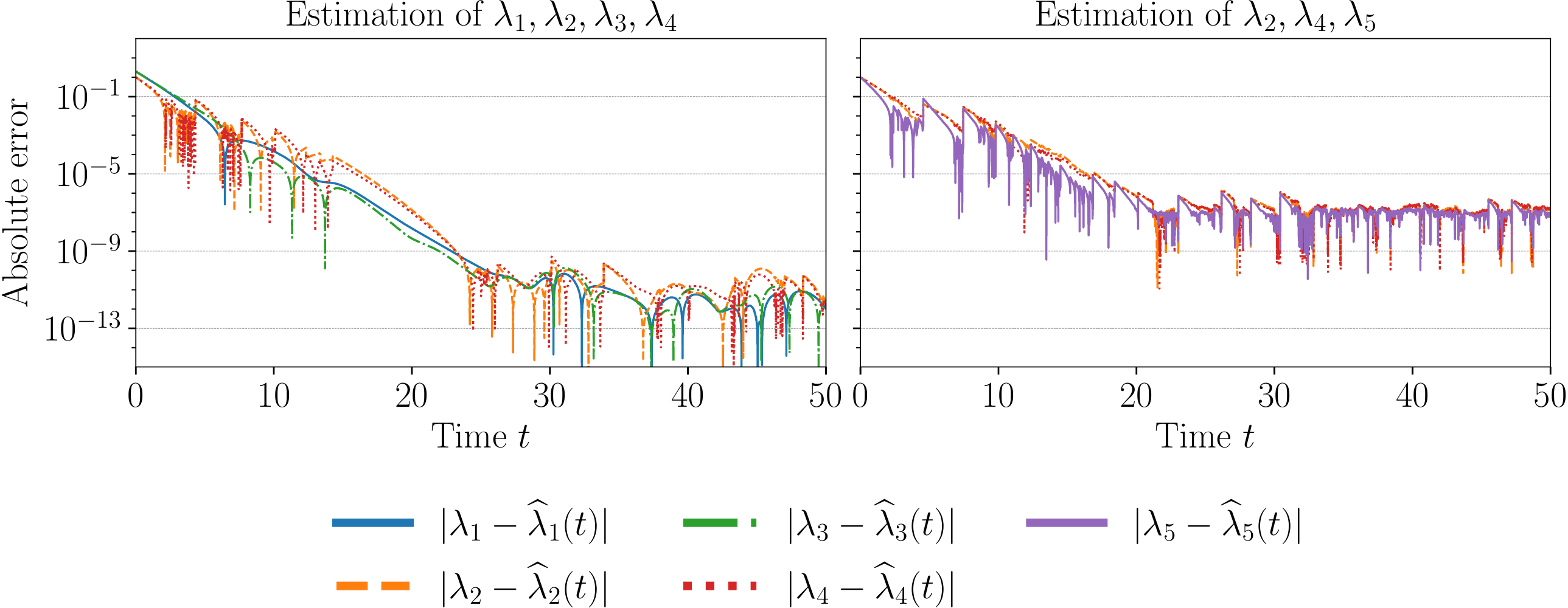}
    \vspace{-0.5cm}
    \caption{Convergence of multi-parameter estimation with respect to time for the four parameters corresponding to linear terms in the generalized KSE \cref{eq:kse_def_generalized} (left) and for the three parameters corresponding to terms in the original KSE, including the coefficient $\lambda_{5}$ on the nonlinear term $uu_x$ (right).}
    \label{fig:convergence_multi}
\end{figure}

As in the case of single-parameter estimation, we observe exponential convergence in time to the true parameters in the multi-parameter case. However, there are additional limitations when multiple parameters are considered. When estimating multiple unknown parameters, the point-in-time estimates must be smoothed as in \cref{eq:relaxation_ode} or the algorithm will fail to converge, i.e. the rapid changes in parameter values introduce non-physical discontinuous features and reactions to the system which overwhelm the continuous data assimilation. Additionally, not every combination of unknown parameters can be learned (at least based on the experiments considered here). In particular, we find that $\lambda_5$ (the coefficient for the nonlinear term) cannot be learned at the same time as $\lambda_1$ or $\lambda_3$---both of which are equal to $0$ in the true KSE. Every other combination of unknown parameters can be successfully recovered. In other words, it appears that this approach can either recover the nonlinear dynamics, or it can perform model discovery, but not both simultaneously.

Estimation of the coefficients of the linear terms $\lambda_1$, $\lambda_2$, $\lambda_3$, and $\lambda_4$ is very efficient and appears to work down to an error of $\mathcal{O}(10^{-11})$, consistent with our previous results when recovering a single parameter. See \cref{fig:convergence_multi} for convergence of all four linear-term coefficients with $\widehat{\lambda}_{5}(t) \equiv \lambda_5 = 1$ known.  Note that the convergence rate is very similar to that observed for the single parameter case (compare with the top right plot of \cref{fig:relaxation}).  It also appears that estimation of the dissipative coefficients ($k=2,4$) behave similar to one another, while the dispersive terms ($k=1,3$) converge in a similar fashion (smoother than the dissipative ones).  This could either be a consequence of the differences between dispersion and dissipation, or it may be simply because $\lambda_1=\lambda_3 = 0$ in the true system.

Complete parameter recovery without the addition of the dispersive $k=1,3$ terms, but including the nonlinear term, is also demonstrated in \cref{fig:convergence_multi}.  Note that the final `converged' error in this case is significantly worse, appearing to level off at $\mathcal{O}(10^{-7})$ for all parameters, even though the convergence rate is very similar (the error stops decreasing at about the same time).  It is not immediately clear why this combination of coefficients does not achieve the same level of accuracy, but it appears that the unknown coefficient on the nonlinearity plays a pivotal role.  Even with this decrease in final accuracy, recovery of the true parameters up to $\mathcal{O}(10^{-7})$ for such a concurrent learning algorithm is significant.

\begin{figure}
    \centering
    \includegraphics[width=\textwidth]{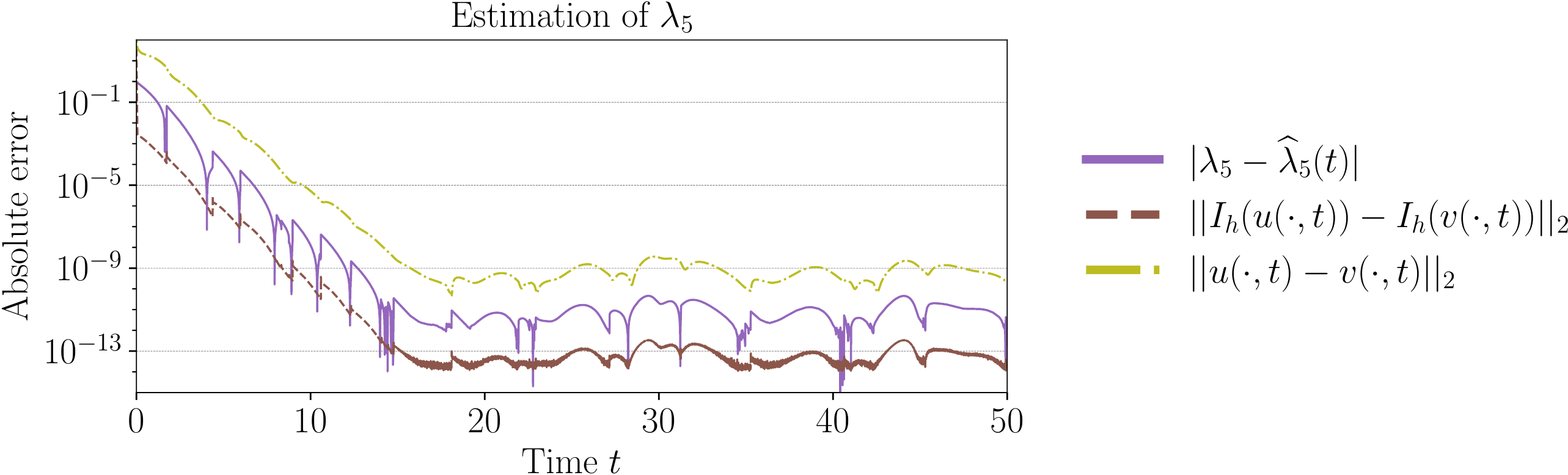}
    \vspace{-0.5cm}
    \caption{Convergence of single-parameter estimation with respect to time for the parameter $\lambda_{5}$ corresponding to the nonlinear term $uu_x$ only in KSE, as well as the state error and the projected state error.}
    \label{fig:convergence_single_nonlin}
\end{figure}

It is worth mentioning that the algorithm experiences no loss of accuracy when $\lambda_5$ is the single unknown parameter, as reported in \cref{fig:convergence_single_nonlin}. A decrease in accuracy only arises when estimating $\lambda_5$ in tandem with other unknown parameters.  This indicates that there is a subtle interplay between recovery of the nonlinear coefficient and recovery of coefficients for the linear terms.

Our experiments indicate that the method appears to be robust to changes in the values of $\alpha$ and/or the initial values of the parameters, i.e., we may make a choice other than $\widehat{\lambda}_k(0) = 2$ with little effect on the convergence of $\widehat{\lambda}_k(t)$ to $\lambda_k$.
On the contrary, the choice of the basis functions $\{e_{i}\}_{i=1}^{N}$ is critical: failing to include the error term $e_1 = I_h(w)$ and instead using basis functions corresponding to the first $n$ Fourier modes leads to a failure to converge. We also implemented an alternative approach that estimated the parameters at each time-step via a least-squares fit using several additional basis functions ($N$ total), as opposed to exactly solving a small linear system. This led to similar convergence to the method presented in detail here, but at slightly increased computational expense. Furthermore, the theoretical justification is less clear for this alternate method, so we leave a more detailed investigation of this and other alternative basis functions to a later study.

\section{Conclusions}
\label{sec:conclusions}

We have developed an algorithm for concurrently learning scalar parameters and assimilating the state in a general setting for evolutionary, dissipative PDEs, and demonstrated that the method works remarkably well on a generalized version of the Kuramoto-Sivashinsky equation.  Numerically, we observe that multiple parameters (coefficients) of KSE can be recovered up to an explainable error (due to the RK time stepping and data assimilation method) so long as these parameters are coefficients on linear terms in the equation.  We have also demonstrated that the coefficient of the nonlinear term can be recovered, but only in concert with dissipative linear terms and then at a nontrivial increase to the final error.  Although there are limitations to the developed approach, including the need for a high temporal resolution, it is shown to work remarkably well on KSE.

Further applications of this algorithm are readily apparent.  Parameter learning or model recovery for physically motivated systems are interesting and well-motivated problems.  Although the derivation in \Cref{sec:algorithm} is not yet rigorous, it does appear that this approach will work for all dissipative systems for which the AOT algorithm has been shown to converge, which represents a very large class of PDEs.  Extensions beyond the simplified setting where the `true' system is run simultaneously as a direct numerical simulation are of greater interest and will be explored at length. For instance, a lofty extension of this work would be to estimate the skin friction of a surface from data collected via wind tunnel experiments.  Such an application of this approach of course entails additional complications that have not been considered here.

Finally, though the derivation provided here is intuitive, it is not yet fully rigorous.  Specifically, we have not provided rigorous justification for the observed convergence of both the parameter learning and the simultaneous synchronization of the state.  The rigorous justification of this algorithm for certain systems of PDEs will not only provide a firm mathematical foundation for this effort, but will also highlight the necessary hypotheses on the hyperparameters (quality of the observations defined by the parameter $h$ for example) that will guide the practical implementation of this method in more practical settings.

\appendix
\section{Numerical Solution to the Generalized KSE}
\label{apdx:solver}
This section provides essential details on the KSE numerical solver.

We consider a generalized version of KSE,
\begin{align*}
    u_{t}
    + \lambda_{1}u_{x}
    + \lambda_{2}u_{xx}
    + \lambda_{3}u_{xxx}
    + \lambda_{4}u_{xxxx} + \lambda_\textrm{5}uu_{x}
    &= 0.
\end{align*}
The classical KSE has $\lambda_{1} = \lambda_{3} = 0$.
Let $\mathcal{F}$ denote the Fourier transform and denote $\widehat{u}(\xi,t) = [\mathcal{F}(u(\cdot,t))](\xi)$.
Recall $[\mathcal{F}(\frac{d^k}{dx^k}\phi)](\xi) = (2\pi i \xi)^{k}[\mathcal{F}(\phi)](\xi)$.
Taking the Fourier transform of both sides of the generalized KSE we arrive at the spectral version of the PDE,
\begin{align*}
    \widehat{u}_{t}
    + \lambda_{1}(2\pi i \xi)\hat{u}
    + \lambda_{2}(2\pi i \xi)^{2}\hat{u}
    + \lambda_{3}(2\pi i \xi)^{3}\hat{u}
    + \lambda_{4}(2\pi i \xi)^{4}\hat{u}
    + \lambda_{5}\mathcal{F}(uu_x)
    = 0.
\end{align*}
To elucidate the numerical implementation of this equation, we let $\mathcal{F}^{-1}$ denote the inverse Fourier transform, i.e. $u(x,t) = \mathcal{F}^{-1}(\widehat{u})$.  This immediately yields the standard pseudo-spectral scheme:
\begin{align*}
    \widehat{u}_t = -\left[\lambda_1 (2\pi \imath \xi) + \lambda_2 (2\pi \imath \xi)^2 + \lambda_3 (2\pi \imath \xi)^3 + \lambda_4 (2\pi \imath\xi)^4\right]\widehat{u}
    \\ \qquad
    - \lambda_5 \mathcal{F}\left(\mathcal{F}^{-1}(\widehat{u})\mathcal{F}^{-1}(2\pi \imath \xi \widehat{u})\right).
\end{align*}
The primary remaining aspects of the numerical scheme are then the choice of a time-stepping discretization of $\widehat{u}_t$ and the finite truncation level of the Fourier transform.

All the simulations performed here are done with an implicit-explicit Runge-Kutta 4-stage method that is formally 4th order, and the Fourier series are all truncated with $N=512$ spatial points.  We also implement a 2/3 dealiasing to avoid unnecessary overflow to the smallest scales.  The choice of domain size and the value of $\lambda_1=1$ is selected for this degree of spatial resolution because even when dealiasing is removed, the power spectrum of the simulation is fully resolved with these parameters, i.e., the magnitude of the smallest-scale terms is near $\epsilon_\textup{machine}$.

For the data assimilation step, the projection operator $I_h(u)$ is defined by zeroing out all modes with wave numbers greater than a cutoff $K$ selected as $K=21$ for all results throughout the paper (except as noted in \cref{sec:projection} and \cref{fig:interpolation}).  All simulations for the true state are started with initial condition
\begin{align*}
    u_0(x)
    &= \sin\left(6\frac{\pi x}{L}\right)
    + 0.1\cos\left(\frac{\pi x}{L}\right)
    - 0.2 \sin \left(3\frac{\pi x}{L}\right)
    \\ & \qquad
    + 0.05 \cos \left(15\frac{\pi x}{L}\right)
    + 0.7\sin\left(18\frac{\pi x}{L}\right)
    - \cos\left(13\frac{\pi x}{L}\right),
\end{align*}
and then run out to non-dimensional time $\tilde{t} = 10$ whereupon the assimilated system is initialized.
\section*{Acknowledgments}
The authors would like to thank E. Carlson, J. Hudson, A. Larios, V. Martinez, and E. Ng for their helpful discussions, as well as the reviewers for their very helpful feedback which led to significant improvements in the manuscript.
JPW was partially supported by the Simons Foundation travel grant 586788.

\bibliographystyle{siamplain}
\bibliography{references}
\end{document}